\documentclass[11pt]{article}

\usepackage{amsmath, amssymb, amsfonts, amsthm, graphicx}

\DeclareMathOperator{\FT}{\mathsf{f\!t}}
\DeclareMathOperator{\bd}{bd}
\DeclareMathOperator{\aff}{aff}
\DeclareMathOperator{\conv}{conv}

\DeclareMathOperator{\interior}{int}
\DeclareMathOperator{\closure}{cl}

\newcommand{\Line}[1]{\overleftrightarrow{\rule{0pt}{1.25ex}#1}}
\newcommand{\Ray}[1]{\overrightarrow{\rule{0pt}{1.25ex}#1}}
\newcommand{\dsegment}[1]{[#1]_d}

\newcommand{\norm}[1]{\|#1\|}
\newcommand{\ipr}[2]{\left\langle #1, #2 \right\rangle}
\newcommand{\un}[1]{\widehat{#1}}

\newcommand{\R}{\mathbf{R}}

\newcommand{\abs}[1]{|#1|}

\newcommand{\myangle}{\sphericalangle}
\newcommand{\vect}[1]{{\boldsymbol{#1}}}
\newcommand{\va}{\vect{a}}
\newcommand{\vb}{\vect{b}}
\newcommand{\vc}{\vect{c}}
\newcommand{\vd}{\vect{d}}

\newcommand{\vp}{\vect{p}}
\newcommand{\vq}{\vect{q}}
\newcommand{\vx}{\vect{x}}
\newcommand{\vy}{\vect{y}}
\newcommand{\vz}{\vect{z}}

\newcommand{\vo}{\vect{o}}

\newcommand{\length}[1]{|#1|}
\newcommand{\opair}[2]{[\![#1,#2]\!]}

\theoremstyle{plain}
\newtheorem{theorem}{Theorem}[section]
\newtheorem{proposition}[theorem]{Proposition}
\newtheorem{lemma}[theorem]{Lemma}
\newtheorem{corollary}[theorem]{Corollary}

\theoremstyle{definition}
\newtheorem{example}[theorem]{Example}
\newtheorem{definition}[theorem]{Definition}

\providecommand{\proof}{\par\textbf{Proof. \quad}}

\begin{document}

\bibliographystyle{amsplain}

\title{The Fermat-Torricelli problem in normed planes and spaces}
\author{Horst Martini \\
        Fakult\"at f\"ur Mathematik \\
        Technische Universit\"at Chemnitz \\
        D-09107 Chemnitz \\ Germany \\
        E-mail: \texttt{martini@mathematik.tu-chemnitz.de} \and 
        Konrad J. Swanepoel\thanks{Research supported by a grant from a cooperation between the Deutsche Forschungsgemeinschaft in Germany and the National Research Foundation in South Africa}\\ 
        Department of Mathematics, Applied Mathematics and Astronomy\\ 
        University of South Africa, \\ P.O.\ Box 392, Unisa 0003 \\ South Africa \\
        E-mail: \texttt{swanekj@unisa.ac.za} \and 
        Gunter Wei\ss \\
        Institut f\"ur Geometrie \\
        Technische Universit\"at Dresden \\
        D-01062 Dresden \\ Germany \\
        E-mail: \texttt{weiss@math.tu-dresden.de}}
\date{}

\maketitle

\section{Introduction}
The famous Fermat-Torricelli problem (in Location Science also called the Steiner-Weber problem) asks for the unique point $\vx$ minimizing the sum of distances to arbitrarily given points $\vx_1,\dots,\vx_n$ in Euclidean $d$-dimensional space $\R^d$.
In the present paper, we will consider the extension of this problem to $d$-dimensional real normed spaces (= Minkowski spaces), where we investigate mainly, but not only, the case $d=2$.

Since in arbitrary Minkowski spaces the solution set (= Fermat-Torricelli locus) is not necessarily a singleton, we study geometric descriptions of this set.
Continuing related investigations given in the papers \cite{CG, DM, DM2, KM}, we present some new geometric results about Fermat-Torricelli loci.
Alongside expositions of known results that are scattered in various sources and proofs of some of them, we make basic observations that have perhaps not been made before, and present many new results, especially in the planar case that is the most important for Location Science.
Our results together can be considered to be a mini-theory of the Fermat-Torricelli problem in Minkowski spaces and especially in Minkowski planes.
We emphasise that it is possible to find substantial results about locational problems valid \emph{for all norms} using a geometric approach, and in fact most of our results are true for all norms.

We now give an overview of the paper.

\textbf{Section~\ref{terminology}.} We introduce and give an overview of basic terminology and technical tools used in Minkowski geometry.

\textbf{Section~\ref{general}.} We first give an overview of general properties of Fermat-Torricelli points and loci (Definitions~\ref{3.1}, \ref{3.2}, \ref{3.3} and Propositions~\ref{prop0} and~\ref{lemma1}).

We then give a simple exposition of the results obtained in \cite{DM} on the characterization of Fermat-Torricelli points in terms of functionals in the dual space (Theorem~\ref{DMchar}) and that the Fermat-Torricelli locus can be obtained as the intersection of certain cones with apices $\vx_1,\dots,\vx_n$ (Definition~\ref{3.7} and Theorem~\ref{geomchar}).
It seems to have been overlooked that this construction is a natural extension of a geometric approach to \emph{$d$-segments} presented in \cite[\S 9]{BMS} (Proposition~\ref{newexample}).
(The notion of $d$-segments was introduced by K. Menger \cite{Menger}, who also gave a historically early investigation in the spirit of Location Science \cite[p.~80]{Menger}.)
We show an application of Theorem~\ref{DMchar} in Corollary~\ref{3.6}.
We then use Theorem~\ref{geomchar} to derive various position criteria for Fermat-Torricelli loci.
In Corollary~\ref{3.9} we describe the shape of Fermat-Torricelli loci in Minkowski planes and give two examples (Examples~\ref{rectilinear} and \ref{reghex}) that will play a role in Section~\ref{specific} in characterizing the $L_1$ and hexagonal norms in the plane (Theorem~\ref{twodim}).

We introduce the new concept of $d$-concurrent $d$-segments (Definition~\ref{3.12}), describe their Fermat-Torricelli loci (Corollary~\ref{dconc}), give examples (Example~\ref{3.14}), and indicate that a result of Cieslik \cite{Cieslik2} follows as a special case (Corollary~\ref{3.15}).
We then introduce the concept of $d$-collinear set (Definition~\ref{3.16}) and characterize Fermat-Torricelli loci of these sets (Corollaries~\ref{dcoll1} and \ref{dcoll2}), thereby generalizing the results in $\R^1$ to general spaces.

In Corollary~\ref{3.19} we describe a general situation of when the Fermat-Torricelli locus is a singleton, and in Theorem~\ref{sc} prove that the Fermat-Torricelli locus is always a singleton exactly when the Minkowski space is strictly convex.

Finally in this section we contrast the situation between the two-dimensional and higher-dimensional cases by citing the result of Wendell and Hurter \cite{WH} that in Minkowski planes the Fermat-Torricelli locus of any set always intersects the convex hull (Theorem~\ref{3.21} and Corollary~\ref{sc2}), and the results of Cieslik \cite{MR89f:90114} and Durier \cite{Durier} that similar properties hold in higher dimensions only in Euclidean space (Theorem~\ref{durier}).
We also sketch the proof of Theorem~\ref{durier}, as its complete proof is scattered over various papers.
In Section~\ref{specific} we refine Theorem~\ref{3.21} (see e.g.\ Theorem~\ref{dc}).

\textbf{Section~\ref{specific}.}
Here we make a closer analysis of the relationship between the Fermat-Torricelli locus and the convex hull of a finite set in Minkowski planes.
We first introduce the notion of a double cluster generalizing the notion of a collinear set with an even number of points (Definition~\ref{4.1} and Example~\ref{doublecluster}) and show that if a Fermat-Torricelli point of a set in a Minkowski plane is outside the convex hull of that set, then the set must be a double cluster (Theorem~\ref{dc}).
It follows that the Fermat-Toricelli locus of a set with an odd number of points is contained in the convex hull of the set (Corollary~\ref{4.4}).

We then introduce the notion of pseudo double cluster generalizing the notion of a set with an even number of points in which all points except possibly one are collinear (Definition~\ref{4.6}) and show that if one of the vertices of the convex hull of a set in a Minkowski plane is also a Fermat-Torricelli point, then the set must be a pseudo double cluster (Theorem~\ref{prop62} and Corollary~\ref{4.8}).
We then give some results on the more subtle situation when there is a Fermat-Torricelli point on the relative interior of an edge of the convex hull (Corollaries~\ref{4.9} and \ref{4.10}, Example~\ref{4.11} and Theorem~\ref{4.12}).
As the final results in Section~\ref{specific} we mention a generalization of Proposition 6.4 of \cite{KM} (Theorem~\ref{4.13}), and characterize the Minkowski planes having parallelograms and affinely regular hexagons as unit balls as those Minkowski planes in which more than two points of a given set can be Fermat-Torricelli points of the set (Theorem~\ref{twodim}).
We give a higher-dimensional generalization of this result (Lemma~\ref{lemma2} and Theorem~\ref{highd}) characterizing $L_1$ spaces.

\textbf{Section~\ref{specific2}.}
Here we conclude our investigation into Fermat-Torricelli loci in Minkowski planes by characterizing absorbing degree two and floating degree three Fermat-Torricelli configurations in terms of special types of angles (absorbing and critical angles, Definitions~\ref{5.1} and \ref{5.2}).
We first give characterizations of these angles (Lemmas~\ref{one} and \ref{two}), which already gives a characterization of absorbing Fermat-Torricelli configurations of degree two.
We then use these results as well as a technical result (Lemma~\ref{threesegment}) to characterize degree three Fermat-Torricelli configurations in terms of critical angles (Theorem~\ref{torithree}).

\section{Terminology of Minkowski spaces}\label{terminology}
A \emph{Minkowski space} is a real finite-dimensional normed space $X$ with \emph{norm} $\norm{\cdot}:X\to\R$ (satisfying $\norm{\vx}\geq 0$, $\norm{\vx}=0$ iff $\vx=\vo$, $\norm{\lambda\vx}=\abs{\lambda}\norm{\vx}$, and most importantly, the triangle inequality $\norm{\vx+\vy}\leq\norm{\vx}+\norm{\vy}$), \emph{unit ball} $B=\{\vx:\norm{\vx}\leq 1\}$ and \emph{unit sphere} (or \emph{unit circle} in the two-dimensional case) $\{\vx:\norm{\vx}= 1\}$.
A \emph{Minkowski plane} is a two-dimensional Minkowski space.
For the facts on Minkowski spaces recalled below, see \cite[Chapters~1~and~3]{Thompson}, for general convex geometry see \cite{Schneider}, and for convex analysis see \cite{Rocka}.

Any centrally symmetric convex body $B$ centred at the origin $\vo$ gives rise to a norm for which $B$ is the unit ball, i.e.,
\[ \norm{\vx} = \inf\{\lambda^{-1} : \lambda\vx\in B\}.\]
By the Mazur-Ulam Theorem, any two Minkowski spaces are isometric iff their unit balls are affinely equivalent, i.e., if there exists a linear mapping from one unit ball onto the other.
A Minkowski space $X$ is \emph{strictly convex} if the unit sphere contains no non-trivial \emph{line segment}
\[\vx\vy=\{\alpha\vx+(1-\alpha)\vy: 0\leq\alpha\leq 1\}, \quad \vx\neq \vy,\]
or, equivalently, if $\norm{\vx+\vy} < \norm{\vx}+\norm{\vy}$ for any linearly independent $\vx, \vy\in X$.
A Minkowski space is \emph{smooth} if each boundary point of the unit ball has a unique supporting hyperplane.

Given a Minkowski space $X$ with norm $\norm{\cdot}$ and unit ball $B$, the dual norm on the dual space $X^\ast$ is defined as $\norm{\phi}=\max_{\norm{\vx}=1}\phi(\vx)$ for any functional $\phi\in X^\ast$.
If we identify $X$ and $X^\ast$ with $d$-dimensional $\R^d$, then the dual unit ball $B^\ast$ is the polar body of $B$:
\[B^\ast = \{\vy: \ipr{\vx}{\vy}\leq 1 \mbox{ for all } \vx\in B\}.\]
A \emph{norming functional} of $\vx\in X$ is a $\phi\in X^\ast$ such that $\norm{\phi}=1$ and $\phi(\vx)=\norm{\vx}$.
The hyperplane $\phi^{-1}(1) = \{\vy\in X : \phi(\vy)=1\}$ is then a hyperplane supporting the unit ball at $\vx$.
By the separation theorem, each $\vx\in X$ has a norming functional.
Thus a Minkowski space is smooth iff each $\vx\neq\vo$ has a unique norming functional.
It is also known that $X^{\ast\ast}$ is isometric to $X$, and $X$ is smooth iff $X^\ast$ is strictly convex.

We use the shorthand notation $\un{\vx}$ for $\frac{1}{\norm{\vx}}\vx$ for any $\vx\neq\vo$, and $\length{\vx\vy}$ for the length $\norm{\vx-\vy}$ of the segment $\vx\vy$.

In some of our proofs we use the subdifferential calculus of convex functions (see \cite[\S 23]{Rocka} for proofs of the discussion below).
A functional $\phi\in X^\ast$ is a {\em subgradient} of a convex function $f: X\to\R$ at $\vx\in X$ if for all $\vz\in X$,
\[ f(\vz)-f(\vx) \geq \phi(\vz-\vx).\]
In particular, $o\in X^\ast$ s a subgradient of $f$ at $\vx$ iff $f$ attains its minimum value at $\vx$.
The {\em subdifferential} of $f$ at $\vx$ is the set $\partial f(\vx)$ of all subgradients of $f$ at $\vx$.
This set is always non-empty, closed and convex.
The following basic property of subdifferentials is important to what follows:
If $f_1,\dots,f_n$ are convex functions on $X$, then
\[ \partial (\sum_{i=1}^n \alpha_i f_i)(\vx) = \sum_{i=1}^n\alpha_i\partial f_i(\vx)\]
for all $\vx\in X$ and $\alpha_1,\dots,\alpha_n\in\R$, where the sum on the right is Minkowski addition of sets in a vector space:
If $A,B\subseteq X$, then $A+B:=\{\va+\vb : \va\in A, \vb\in B\}$.
The proof, to be found in \cite{Rocka}, uses the separation theorem.
It is easily seen that the subdifferential of the norm of $X$ at $\vx$ is the following:
\begin{lemma}\label{normsubdiff}
$\partial\norm{\vo} = B_{X^\ast}$ \textup{(}i.e. the unit ball of $X^\ast$\textup{)}.
If $\vx\neq \vo$, then $\partial\norm{\vx} = \{\phi\in X^\ast: \norm{\phi} = 1, \phi(\vx)=\norm{\vx} \}$, \textup{(}i.e. the set of norming functionals of $\vx$\textup{)}.
\qed
\end{lemma}
Thus, if $\vx\neq \vo$, $\partial\norm{\vx}$ is the exposed face of the unit ball in $X$ defined by the hyperplane $\{\phi\in X^\ast: \phi(\vx) = 1\}$.
Recall that a (proper) \emph{exposed face} of a convex body $B$ is an intersection of $B$ with some supporting hyperplane (see e.g.\ \cite{Schneider}).

We conclude with some more geometric terms.
The \emph{ray} with origin $\va$ passing through $\vb$ is denoted by $\Ray{\va\vb}$.
An \emph{angle} $\myangle \vx\vy\vect{z}$ in a Minkowski plane is the convex cone bounded by two rays $\Ray{\vy\vx}$ and $\Ray{\vy\vect{z}}$ emanating from the same point $\vy$.
(We allow half planes, i.e.\ $180^\circ$ angles --- in this case we take the half plane on the left if we pass from $\vx$ to $\vz$.)
We denote the \emph{$d$-segment} from $\vx$ to $\vy$ by
\[\dsegment{\vx\vy} := \{\vz\in X:\length{\vx\vz} + \length{\vz\vy} = \length{\vx\vy}\}.\]
A \emph{metric ray} is a subset of $X$ that is isometric to $[0,\infty)$, and a \emph{metric line} is a subset of $X$ isometric to $\R^1$.

We denote the interior, closure, boundary, convex hull and affine hull of a subset $A$ of a Minkowski space by $\interior A$, $\closure A$, $\bd A$, $\conv A$ and $\aff A$, respectively.

\section{Fermat-Torricelli points and loci: General properties}\label{general}
\begin{definition}\label{3.1}
We call a point $\vx_0$ a \emph{Fermat-Torricelli point} (or FT point) of distinct points $\vx_1, \dots, \vx_n$ in a Minkowski space if $\vx=\vx_0$ minimizes $\vx\mapsto\sum_{i=1}^n\length{\vx\vx_i}$.
\end{definition}
See \cite{KM} and \cite[Chapter 2]{MR2000c:90002} for a discussion of FT points in Euclidean spaces; related investigations in Minkowski spaces are \cite{CG, Durier, DM, S6}; see also \cite{Cieslik2}.
In the Facilities Location literature (cf.\ \cite{MR1358610}) these points are also called Fermat-Weber or Steiner-Weber points.
\begin{definition}\label{3.2}
A (star) \emph{configuration} (of degree $n$) in a Minkowski space is a set of segments $\{\vx\vx_i: i=1, \dots,n\}$ emanating from the same point $\vx$, with $\vx_i\neq\vx$ for all $i$.
A configuration $\{\vx\vx_i\}$ is \emph{pointed} if there is a hyperplane $H$ through $\vx$ such that the interior of each segment $\vx\vx_i$ is in the same open half space bounded by $H$.
A \emph{floating Fermat-Torricelli configuration} (or floating FT configuration) is a configuration $\{\vx_0\vx_i:i=1, \dots,n\}$ such that $\vx_0$ is an FT point of $\{\vx_1, \dots, \vx_n\}$, and an \emph{absorbing Fermat-Torricelli configuration} (or absorbing FT configuration) is a configuration $\{\vx_0\vx_i:i=1, \dots,n\}$ such that $\vx_0$ is an FT point of $\{\vx_0, \vx_1, \dots, \vx_n\}$.
\end{definition}
We first make the following simple observations.
\begin{proposition}\label{prop0}
In any Minkowski space,
\begin{enumerate}
\item if $\{\vx_0\vx_i\}$ is a floating FT configuration, then it is also an absorbing FT configuration.
\item if $\{\vx_0\vx_i\}$ is an FT configuration, then so is $\{\vx_0\vy_i\}$ for any $\vy_i\in\Ray{\vx_0\vx_i}$, $\vy_i\neq \vx_0$.
\end{enumerate}
\end{proposition}
\proof
Firstly, if $\vx=\vx_0$ minimizes $\vx\mapsto\sum_{i=1}^n\length{\vx\vx_i}$, then $\vx=\vx_0$ also minimizes $\vx\mapsto\sum_{i=0}^n\length{\vx\vx_i}$, since for any $\vx\in X$ we have $\sum_{i=0}^n\length{\vx_0\vx_i}=\sum_{i=1}^n\length{\vx_0\vx_i}\leq\sum_{i=1}^n\length{\vx\vx_i}\leq\sum_{i=0}^n\length{\vx\vx_i}$.

Secondly, suppose $\{\vx_0\vx_i : i=1,\dots,n\}$ is a floating FT configuration.
(The case of an absorbing FT configuration is similar.)
Without loss of generality $\vx_0=\vo$.
Thus $\vx=\vo$ minimizes $\vx\mapsto\sum_{i=1}^n\length{\vx\vx_i}$.
Let $\vy_i\neq\vo$ be on the ray $\Ray{\vo\vx_i}$ for each $i=1,\dots,n$, say $\vy_i=\lambda_i\vx_i$ with $\lambda_i>0$.
Clearly $\{\vo\vx_i\}$ is an FT configuration iff $\{\vo\vx_i'\}$ is an FT configuration where $\vx_i'=\lambda\vx_i$, for any $\lambda>0$, i.e., we may scale an FT configuration.
Thus we may assume without loss of generality that each $\lambda_i\leq 1$ by making the original FT configuration sufficiently large.
Then for any $\vx\in X$ we have
\begin{eqnarray*}
\sum_{i=1}^n\length{\vo\vy_i} & = & \sum_{i=1}^n(\length{\vo\vx_i}-\length{\vx_i\vy_i}) \\
& \leq & \sum_{i=1}^n(\length{\vx\vx_i}-\length{\vx_i\vy_i}) \quad \text{ (since $\vo$ is an FT point)} \\
& \leq & \sum_{i=1}^n(\length{\vx\vy_i} \quad \text{ (by the triangle inequality).}
\end{eqnarray*}
Thus $\vo$ is an FT point of $\{\vy_i\}$.
\qed

In contrast to the case of non-collinear points in Euclidean space, in general a set of points can have more than one FT point.
\begin{definition}\label{3.3}
The \emph{Fermat-Torricelli locus} (or FT locus) of $\vx_1, \dots, \vx_n$ is the set $\FT(\vx_1, \dots, \vx_n)$ of all FT points of $\vx_1, \dots, \vx_n$.
\end{definition}
There are also corresponding definitions for weighted points, but we only consider the unweighted case.
Note that it immediately follows from the triangle inequality that $\FT(\vx,\vy)=\dsegment{\vx\vy}$.
Chakerian and Ghandehari \cite{CG} gave an extensive analysis of FT points in the floating case, where $X$ is a smooth and strictly convex Minkowski space.
They derive characterizations in terms of ``special polytopes'', i.e., in the terminology of \cite{KM}, polyhedral arrangements with the Viviani-Steiner property.
In our discussion we do not make in general any special assumptions such as smoothness or strict convexity, nor do we restrict our attention exclusively to the floating case.
In general we can say the following of the FT locus \cite{MR89f:90114}.

\begin{proposition}\label{lemma1}
The FT locus of any finite set is always non-empty, compact and convex.
\end{proposition}
\proof
The following is a standard argument, adapted from the Euclidean case.
If $A=\{\vx_1,\dots,\vx_n\}$, then $\FT(A)$ is the set of all minima of the function
\[f(\vx):=\sum_{i=1}^n\length{\vx\vx_i}.\]
By the triangle inequality,
\[f(\vx)\geq\sum_{i=1}^n(\norm{\vx}-\norm{\vx_i})>\sum_{i=1}^n\norm{\vx_i}=f(\vo)\]
for any $\vx$ with $\norm{\vx}>2\sum_{i=1}^n\norm{\vx_i}$.
Thus $\FT(A)$ is contained in the closed ball $\norm{\vx}\leq 2\sum_{i=1}^n\norm{\vx_i}$, and by compactness $\FT(A)$ is non-empty and compact.
That $\FT(A)$ is convex follows from the convexity of the function $f$.
\qed

Durier and Michelot \cite{DM} have given the following characterization of FT points, which extends the classical characterization in the case of Euclidean spaces.

\begin{theorem}[\cite{DM}]\label{DMchar}
Let $\vx_0, \vx_1, \dots, \vx_n$ be points in a Minkowski space.
\begin{enumerate}
\item If $\vx_0\neq\vx_1, \dots, \vx_n$, then $\{\vx_0\vx_i:i=1, \dots,n\}$ is a floating FT configuration iff each $\vx_i-\vx_0$ has a norming functional $\phi_i$ such that $\sum_{i=1}^n\phi_i=o$.
\item If $\vx_0=\vx_j$ for some $j=1, \dots,n$, then $\{\vx_0\vx_i:i=1, \dots,n, i\neq j\}$ is an absorbing FT configuration iff each $\vx_i-\vx_0$ $(i\neq j)$ has a norming functional $\phi_i$ such that
\[\Biggl\|\sum_{\substack{i\neq j\\ i=1}}^n\phi_i\Biggr\|\leq 1.\]
\end{enumerate}
\end{theorem}

\proof
Let $A=\{\vx_0,\dots,\vx_n\}$.
We use the subdfferential calculus.
The point $\vp\in\FT(A)$ iff $\vp$ minimizes the convex function
\[ f(\vx)=\sum_{i=1}^n\length{\vx\vx_i}, \]
iff $\vo\in\partial f(\vx)=\partial\sum_{i=1}^n\length{\vx\vx_i}=\sum_{i=1}^n\partial\length{\vx\vx_i}$.
This is equivalent to the conditions stated, since, letting $g(\vx)=\length{\vx\vx_i}$, we have by Lemma~\ref{normsubdiff} that
\[ \partial g(\vx) = \left\{ \begin{array}{l@{\mbox{ if }}l}
\{\phi: \phi \mbox{ is a norming functional of } \vx-\vx_i\} & \vx\neq \vx_i \\
\{\phi: \norm{\phi}\leq 1\} & \vx=\vx_i.\end{array}\right. \]

Sufficiency can also be shown directly as follows:
If the purported FT point $\vp\not\in A$ then for any $\vx\in X$,
\begin{eqnarray*}
\sum_{i=1}^n\length{\vx_i\vp} & = & \sum_{i=1}^n\phi_i(\vx_i-\vp) \\ 
& = & \sum_{i=1}^n\phi_i(\vx_i-\vx)+\sum_{i=1}^n\phi_i(\vx-\vp) \\ 
& = & \sum_{i=1}^n\phi_i(\vx_i-\vx)+(\sum_{i=1}^n\phi_i)(\vx-\vp) \\
& = & \sum_{i=1}^n\phi_i(\vx_i-\vx) \\
& \leq & \sum_{i=1}^n\length{\vx_i\vx},
\end{eqnarray*}
while if $\vp=\vx_j\in A$, then for any $\vx\in X$,
\begin{eqnarray*}
\sum_{i\neq j}\length{\vx_i\vp}  & = & \sum_{i\neq j}\phi_i(\vx_i-\vp) \\
& = & \sum_{i\neq j}\phi_i(\vx_i-\vx)+\sum_{i\neq j}\phi_i(\vx-\vp) \\
& = & \sum_{i\neq j}\phi_i(\vx_i-\vx)+(\sum_{i\neq j}\phi_i)(\vx-\vp) \\
& \leq & \sum_{i\neq j}\length{\vx_i\vx}+\Bigl\|\sum_{i\neq j}\phi_i\Bigr\|\,\length{\vx\vp} \\
& \leq & \sum_{i=1}^n\length{\vx_i\vx}.
\end{eqnarray*}
\qed

The above type of calculation is useful to analyze the situation where $\FT(A)$ has more than one point (see the proofs of Theorem~\ref{sc} and Lemma~\ref{lemma2}).
Theorem~\ref{prop0} follows immediately from the above characterization, as well as the following observation.
\begin{corollary}\label{3.6}
In any Minkowski space, if $\{\vo\vx_1,\dots,\vo\vx_n\}$ is a floating FT configuration, then $\{\vo\vx_1,\dots,\vo\vx_{n-1}\}$ is an absorbing FT configuration.
\end{corollary}

The following geometric description of $\FT(A)$ (Theorem~\ref{geomchar}), due to Durier and Michelot \cite{DM}, also follows from Theorem~\ref{DMchar}.
\begin{definition}\label{3.7}
Given a unit functional $\phi\in X^*$ and a point $\vx\in X$, define the cone $C(\vx,\phi)=\vx-\{\va:\phi(\va)=\norm{\va}\}$, i.e., $C(\vx,\phi)$ is the translate by $\vx$ of the union of the rays from the origin through the exposed face $\phi^{-1}(-1)\cap B$ of the unit ball $B$ of $X$.
\end{definition}
We remark that any metric ray of $X$ with origin $\vx$ is contained in $C(\vx,\phi)$ for some unit functional $\phi$.

Note that it follows from Proposition~\ref{lemma1} that if $\FT(A)$ consists of more than one point then it contains a point not in $A$, since then $\FT(A)$ is infinite, but $A$ is finite.
In the following theorem we need an FT point not in $A$ in order to obtain a geometric description of $\FT(A)$.
To apply this theorem we therefore first have to find such an FT point by some other means.
\begin{theorem}[\cite{DM}]\label{geomchar}
In any Minkowski space $X$ with a finite given subset $A$, suppose we are given $\vp\in\FT(A)\setminus A$.
Let $\phi_i$ be a norming functional of $\vx_i-\vp$ for each $\vx_i\in A$, such that $\sum_i\phi_i=o$.
Then $\FT(A)=\bigcap_{i=1}^n C(\vx_i,\phi_i)$.
\end{theorem}
\proof
Note that by Definition~\ref{3.7} $\vx\in\bigcap_{i=1}^n C(\vx_i,\phi_i)$ iff $\phi_i(\vx_i-\vx)=\norm{\vx_i-\vx}$.
Thus if $\vx\not\in A$ we have that $\vx\in\bigcap_{i=1}^n C(\vx_i,\phi_i)$ iff for each $i=1,\dots,n$, $\phi_i$ is a norming functional of $\vx_i-\vx$, iff $\vx\in\FT(A)$ (by Theorem~\ref{DMchar} and $\sum_i\phi_i=o$).
It follows that $\FT(A)\setminus A = \bigcap_{i=1}^n C(\vx_i,\phi_i) \setminus A$.

If on the other hand $\vx=\vx_j$ for some $j$, then $\vx\neq\vx_i$ for all $i\neq j$, and $\vx\in\bigcap_{i=1}^n C(\vx_i,\phi_i)$ implies that for all $i\neq j$, $\phi_i$ is a norming functional of $\vx_i-\vx$, which implies that $\vx\in\FT(A)$ (by Theorem~\ref{DMchar} and $\norm{\sum_{i\neq j}\phi_i}=\norm{-\phi_j}=1$).
Thus $\bigcap_{i=1}^n C(\vx_i,\phi_i)\subseteq\FT(A)$.

It remains to show that $A\cap\FT(A)\subseteq\bigcap_{i=1}^n C(\vx_i,\phi_i)$.
Since $\FT(A)\setminus A$ is not empty, and $\FT(A)$ is convex, hence connected, we have that any $\vx_i\in A\cap\FT(A)$ is not an isolated point of $\FT(A)$.
Thus we have $\vx_i\in\closure(\FT(A)\setminus A)=\closure(\bigcap_{i=1}^n C(\vx_i,\phi_i) \setminus A)\subseteq\bigcap_{i=1}^n C(\vx_i,\phi_i)$.
\qed

In the special case where $A$ consists of two points $\vx$ and $\vy$, $\FT(A)$ is the $d$-segment $\dsegment{\vx\vy}$, and from the above theorem can be found the description of $d$-segments obtained in \cite[Theorem~9.6]{BMS}.
We here demonstrate the planar case.
Since we may make a translation we assume without loss of generality in the following proposition that $\vy=-\vx$.

\begin{figure}
\begin{center}
\includegraphics[scale = 0.9]{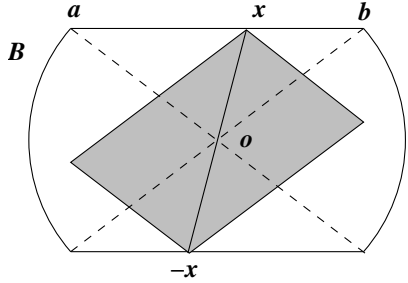}
\end{center}
\caption{A $d$-segment in a Minkowski plane}\label{dsegmentfig}
\end{figure}

\begin{proposition}\label{newexample}
In a Minkowski plane for any $\vx\neq\vo$, we have that $\FT(\vx,-\vx)=\dsegment{\vx\vy}$ is the segment $\vx(-\vx)$ whenever $\vx$ is not in the relative interior of a segment on the boundary of the unit ball, while if $\vx$ is in the relative interior of the maximal segment $\va\vb$ on the boundary of the unit ball, then $\FT(\vx,-\vx)$ is the \textup{(}unique\textup{)} parallelogram with sides parallel to $\Line{\vo\va}$ and $\Line{\vo\vb}$ and which has $\vx$ and $-\vx$ as opposite sides.
See Figure~\ref{dsegmentfig}.
\end{proposition}
\proof
Take any norming functional $\phi$ of $\vx$.
Then $\phi^{-1}(1)\cap B$ contains $\vx$ and is either the singleton $\{\vx\}$ or a segment $\va\vb$.
In the case of a singleton clearly $C(\vx,\phi)$ is the ray $\Ray{\vx(-\vx)}$.
Since $-\phi$ is then a norming functional of $-\vx$ and $C(-\vx,-\phi)=\Ray{(-\vx)\vx}$, if we apply Theorem~\ref{geomchar}, we obtain $\FT(\vx,-\vx)=\vx(-\vx)$.

If $\phi^{-1}(1)\cap B=\va\vb$, then $C(\vx,\phi)$ is the angle with vertex $\vx$ bounded by two rays with origin $\vx$ in the directions of $-\va$ and $-\vb$, respectively.
Applying Theorem~\ref{geomchar} we obtain that $\FT(\vx,-\vx)$ is the intersection of these two angles, which is the parallelogram described in the statement of the proposition.

Note that if $\vx$ is not in the relative interior of $\va\vb$, i.e., if $\vx=\va$ or $\vx=\vb$, then the parallelogram degenerates to the segment $\vx(-\vx)$.
\qed

As seen in the above proof, in a Minkowski plane the cones $C(\vx_i,\phi)$ are always either rays or angles, and in the light of Proposition~\ref{lemma1} we obtain the following
\begin{corollary}\label{3.9}
In a Minkowski plane $X$, the FT locus of a finite set of points is always a convex polygon, that may degenerate to a segment or a point.
\end{corollary}

We here give two examples in detail of how Theorems~\ref{DMchar} and \ref{geomchar} can be applied to find FT loci.
Later it will be seen that these examples are unique in a certain sense (Theorem~\ref{twodim}).
\begin{example}\label{rectilinear}
Let the unit ball of the Minkowski plane $X$ be the parallelogram with vertices $\{\pm\vx,\pm\vy\}$, where $\vx$ and $\vy$ are any two linearly independent vectors.
If $\vx$ and $\vy$ form the standard basis of $\R^2$, then we obtain the $L_1$ plane or Manhattan plane.
In fact, a Minkowski plane is isometric to the $L_1$ plane iff its unit ball is a parallelogram.
We let $A=\{\pm\frac12(\vx+\vy), \pm\frac12(\vx-\vy)\}$ (see Figure~\ref{fig1}).
We now use Theorem~\ref{geomchar} to find $FT(A)$.
\begin{figure}
\begin{center}
\includegraphics[scale = 0.9]{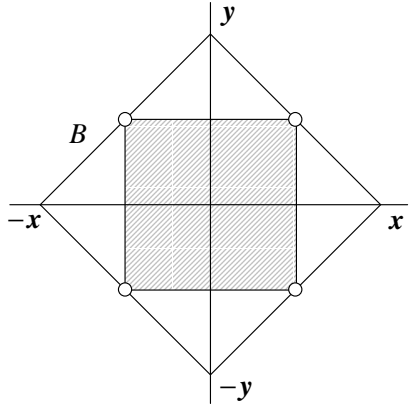}
\end{center}
\caption{$\FT(A)=\conv(A)$ is possible in the rectilinear norm}\label{fig1}
\end{figure}
We let $\phi_1$ be the norming functional of $\frac12(\vx+\vy)$, i.e., the (unique) functional in $X^*$ of norm $1$ for which $\phi_1^{-1}(1)=\Line{\vx\vy}$.
Similarly, we let $\phi_2$ be the norming functional of $\frac12(\vx-\vy)$, i.e., the functional of norm $1$ for which $\phi_2^{-1}(1)=\Line{\vx(-\vy)}$.
Then the norming functionals of $-\frac12(\vx+\vy)$ and $-\frac12(\vx-\vy)$ are $-\phi_1$ and $-\phi_2$, respectively.
By Theorem~\ref{DMchar} we then have that $\vo$ is an FT point of $A$ (since the sum of the norming functionals is $o$).
By Theorem~\ref{geomchar}  we have
\begin{eqnarray*}
\FT(A) & = & C(\tfrac12(\vx+\vy),\phi_1)\cap C(\tfrac12(\vx-\vy),\phi_2) \\
& & \;\cap \,C(-\tfrac12(\vx+\vy),-\phi_1)\cap C(-\tfrac12(\vx-\vy),-\phi_2).
\end{eqnarray*}
The union of the rays from the origin through the exposed face $\phi_1^{-1}(-1)\cap B$ of the unit ball is the whole third quadrant $\myangle(-\vx)\vo(-\vy)$.
Thus $C(\frac12(\vx+\vy),\phi_1)$ is the translate of this quadrant by $\frac12(\vx+\vy)$, i.e., the angle $\myangle\frac12(-\vx+\vy)\frac12(\vx+\vy)\frac12(\vx-\vy)$.
The other cones are similarly found, and their intersection is exactly the parallelogram $\conv A$, which is the shaded part in Figure~\ref{fig1}.
\end{example}

\begin{example}\label{reghex}
Let the unit ball of the Minkowski plane $X$ be an affine regular hexagon $B$, i.e., $B$ is the image of a regular hexagon with centre $\vo$ under an invertible linear mapping.
If we let $\vx$ and $\vy$ be two consecutive vertices of $B$ then $B=\conv\{\pm\vx,\pm\vy,\pm(\vx-\vy)\}$.
See Figure~\ref{fig2}.
We now use Theorem~\ref{geomchar} to find $\FT(\vo,\vx,\vy)$.
\begin{figure}
\begin{center}
\includegraphics[scale = 0.9]{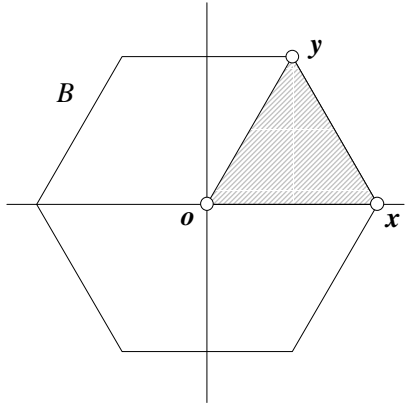}
\end{center}
\caption{$\FT(A)=\conv(A)$ is possible in the regular hexagonal norm}\label{fig2}
\end{figure}
Let $\vp=\frac13(\vo+\vx+\vy)$, i.e., the centroid of the triangle $\triangle\vo\vx\vy$.
As in Example~\ref{rectilinear}, if we let $\phi_1,\phi_2,\phi_3$ be the (unique) norming functionals of $\vo-\vp, \vx-\vp, \vy-\vp$, respectively, then $\phi_1+\phi_2+\phi_3=o$.
Thus $\vp$ is an FT point of $\{\vo,\vx,\vy\}$.
As before, we have $C(\vo,\phi_1)=\myangle\vx\vo\vy$, $C(\vx,\phi_2)=\myangle\vo\vx\vy$, $C(\vy,\phi_3)=\myangle\vo\vy\vx$.
By Theorem~\ref{geomchar}, $\FT(\vo,\vx,\vy)=\myangle\vx\vo\vy\cap\myangle\vo\vx\vy\cap\myangle\vo\vy\vx=\conv\{\vo,\vx,\vy\}$,
i.e., the FT locus is the triangle $\triangle\vo\vx\vy$.
\end{example}

\begin{definition}\label{3.12}
The $d$-segments $\dsegment{\va_i\vb_i}$ are \emph{$d$-concurrent} if their intersection is non-empty.
\end{definition}
\begin{corollary}\label{dconc}
If $A=\{\vx_1,\dots,\vx_{2k}\}$ can be matched up to form $k$ $d$-segments $\dsegment{\vx_i\vx_{k+i}}$, $i=1,\dots,k$, that are $d$-concurrent, then $\FT(A)=\bigcap_{i=1}^{k}\dsegment{\vx_i\vx_{k+i}}$.
\end{corollary}
\proof
Let $\vp\in\bigcap_{i=1}^{k}\dsegment{\vx_i\vx_{k+i}}\setminus A$.
Since $\dsegment{\vx_i\vx_{k+i}}=\FT(\vx_i,\vx_{k+i})$, we have that $\{\vp\vx_i, \vp\vx_{k+1}\}$ is a floating FT configuration for each $i$.
By Theorem~\ref{DMchar} there is a norming functional $\phi_i$ of $\vx_i-\vp$ and $\phi_{k+i}$ of $\vx_{k+i}-\vp$ such that $\phi_i+\phi_{k+i}=\vo$.
By Theorem~\ref{geomchar} we have $\FT(\vx_i\vx_{k+i})=C(\vx_i,\phi_i)\cap C(\vx_{k+i},\phi_{k+i})$.
Thus $\sum_{i=1}^{2k}\phi_i=o$, and again by Theorem~\ref{DMchar}, $\{\vp\vx_1,\dots,\vp\vx_{2k}\}$ is a floating FT configuration.
We now apply Theorem~\ref{geomchar} again to obtain $\FT(A)=\bigcap_{i=1}^{2k}C(\vx_i,\phi_k)=\bigcap_{i=1}^k\FT(\vx_i\vx_{k+i})$.
\qed

\begin{example}\label{3.14}
In Figure~\ref{fig4} we first apply Proposition~\ref{newexample} to obtain $\dsegment{\vx_1\vx_3}$ and $\dsegment{\vx_2\vx_4}$.
We obtain that $\dsegment{\vx_1\vx_3}$ is the usual segment $\vx_1\vx_3$, since the exposed faces of the unit ball $B$ containing $\vx_1$ and $\vx_3$, respectively, are both singletons ($\vx_1$ and $\vx_3$ are not in the interiors of segments on the boundary of $B$).
Also, $\dsegment{\vx_2\vx_4}$ is the parallelogram with opposite sides $\vx_2$ and $\vx_4$ and sides parallel to the vectors from the origin to the endpoints of the segment on the boundary of $B$ containing $\vx_4$.
By Corollary~\ref{3.12}, $\FT(\vx_1,\vx_2,\vx_3,\vx_4)=\dsegment{\vx_1\vx_3}\cap\dsegment{\vx_2\vx_4}$ (since the intersection is non-empty), which is the diagonal of the parallelogram $\dsegment{\vx_2\vx_4}$ indicated as a bold segment in Figure~\ref{fig4}.
\begin{figure}
\begin{minipage}{6.6cm}
\begin{center}
\includegraphics[scale = 0.85]{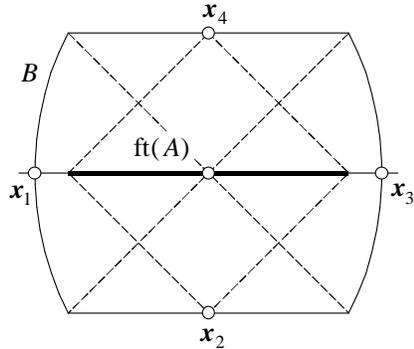}
\end{center}
\caption{$d$-concurrent $d$-segments}\label{fig4}
\end{minipage}
\hfill
\begin{minipage}{6cm}
\begin{center}
\includegraphics[scale = 0.99]{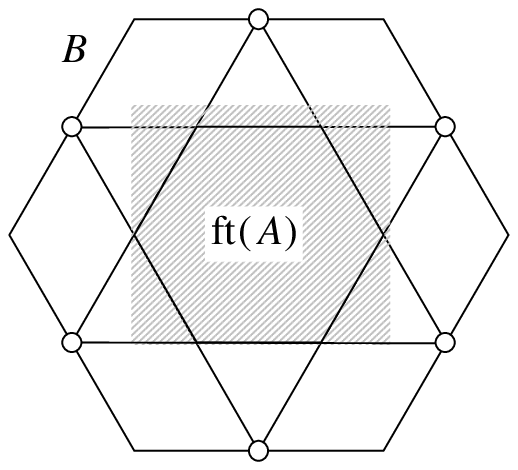}
\end{center}
\caption{$d$-concurrent $d$-segments}\label{fig3}
\end{minipage}
\end{figure}

In Figure~\ref{fig3} the Minkowski plane has an affine regular hexagon as unit ball $B$.
The set $A$ consists of the midpoints of the edges of $B$.
If we now apply Proposition~\ref{newexample} to pairs of points on opposite edges, we obtain that the $d$-segments of these pairs of points has non-empty intersection, which is the shaded hexagon in Figure~\ref{fig3}.
By Corollary~\ref{3.12}, this hexagon is $\FT(A)$.
\end{example}

The following result of Cieslik \cite{Cieslik2} generalizes the Euclidean case \cite{KM}, and follows from Corollary~\ref{dconc}.

\begin{corollary}[{\cite[Chapter~3]{Cieslik2}}]\label{3.15}
Let $\va\vb\vc\vd$ be a convex quadrilateral in a Minkowski plane.
Then the intersection of the diagonals $\va\vc\cap\vb\vd$ is an FT point of $\{\va,\vb,\vc,\vd\}$.
\end{corollary}

The next two corollaries generalize the standard results on FT points in $\R^1$.
\begin{definition}\label{3.16}
A set in a Minkowski space is \emph{$d$-collinear} if it is contained in a metric line.
\end{definition}
Since two finite isometric subsets of $\R^1$ differ by a translation and possibly a reflection, the ordering of a $d$-collinear set is essentially unique.
\begin{corollary}\label{dcoll1}
If $A=\{\vx_1,\dots,\vx_{2k}\}$ is a $d$-collinear set of even size in its natural order, then
\[\FT(A)=\bigcap_{i=1}^k\dsegment{\vx_i\vx_{2k-i}} = \dsegment{\vx_k\vx_{k+1}}.\]
\end{corollary}
\proof
The first equation follows from the observation that $\vx_j\in\dsegment{\vx_i\vx_k}$ for all $1\leq i\leq j\leq k\leq n$, which holds since $\length{\vx_i\vx_k}=\length{\vx_i\vx_j}+\length{\vx_j\vx_k}$, since $\{\vx_i,\vx_j,\vx_k\}$ is isometric to a subset $\{r_1,r_2,r_3\}$ of $\R^1$ with $r_1\leq r_2\leq r_3$.
Thus the elements of $A$ can be matched up to form $d$-concurrent $d$-segments $\dsegment{\vx_i\vx_{2k-i+1}}$, and Corollary~\ref{3.12} applies.
The second equation follows from the fact that $\dsegment{\vx_k\vx_{k+1}}\subseteq\dsegment{\vx_i\vx_{2k-i+1}}$ for all $i=1,\dots,k$, which in turn is proved as follows:
Let $\vx\in\dsegment{\vx_k\vx_{k+1}}$.
Thus $\length{\vx_k\vx}+\length{\vx\vx_{k+1}}=\length{\vx_k\vx_{k+1}}$.
Then
\begin{eqnarray*}
\length{\vx_i\vx} + \length{\vx\vx_{2k-i+1}} & \leq & \length{\vx_i\vx_k}+\length{\vx_k\vx}+\length{\vx\vx_{k+1}}+\length{\vx_{k+1}\vx_{2k-i+1}} \\
& & \qquad \text{ (by the triangle inequality)} \\
& = & \length{\vx_i\vx_k}+\length{\vx_k\vx_{k+1}}+\length{\vx_{k+1}\vx_{2k-i+1}} \\
& = & \length{\vx_i\vx_{2k-i+1}} \\
& & \qquad \text{ (since $\vx_i,\vx_k,\vx_{k+1},\vx_{2k-i+1}$ are on a metric line).}
\end{eqnarray*}
By the triangle inequality we then have $\length{\vx_i\vx} + \length{\vx\vx_{2k-i+1}}=\length{\vx_i\vx_{2k-i+1}}$, hence $\vx\in\dsegment{\vx_i\vx_{2k-i+1}}$.
\qed

\begin{corollary}\label{dcoll2}
If $A=\{\vx_1,\dots,\vx_{2k+1}\}$ is a $d$-collinear set of odd size in its natural order, then $\FT(A)=\{\vx_k\}$.
\end{corollary}
\proof
By Corollary~\ref{dcoll1}, $\vx_k$ is an FT point of $A\setminus\{\vx_k\}$.
Therefore, for any $\vp\in X$,
\begin{eqnarray*}
\sum_{i=1}^{2k+1}\length{\vx_i\vp} & = & \length{\vx_k\vp} + \sum_{\substack{i=1\\ i\neq k}}^{2k+1}\length{\vx_i\vp} \\
& \geq & \length{\vx_k\vp} + \sum_{i=1}^{2k+1}\length{\vx_i\vx_k} \\
& > & \sum_{i=1}^{2k+1}\length{\vx_i\vx_k} \mbox{ unless } \vp=\vx_k.
\end{eqnarray*}
\qed

The calculation in the above proof also gives the following simple

\begin{corollary}\label{3.19}
If $\vp\in\FT(A)$ and $\vp\not\in A$, then $\FT(A\cup\{\vp\})=\{\vp\}$.
\end{corollary}

It is well-known that there is a unique FT point for any non-collinear set in Euclidean space.
The essential property of Euclidean space that ensures uniqueness is its strict convexity.

\begin{theorem}\label{sc}
A Minkowski space $X$ is strictly convex  iff $\FT(A)$ is a singleton for all non-collinear subsets $A$.
\end{theorem}
\proof
If $X$ is not strictly convex, and if we let $\va,\vb$ be two points in the relative interior of some segment on the boundary of the unit ball, then $\FT(\va,-\va,\vb,-\vb)$ is not a singleton by Corollary~\ref{dcoll1}.

Conversely, suppose that $\vp$ and $\vq$ are distinct FT points of a finite set $A$.
Then the segment $\vp\vq\subseteq\FT(A)$, by Proposition~\ref{lemma1}.
Thus we may assume that $\vp,\vq\not\in A$, since $A$ is finite.
By Theorem~\ref{DMchar}, there exist norming functionals $\phi_i$ of $\vx_i-\vp$ for each $\vx_i\in A$ such that $\sum_i\phi_i=o$.
Thus,
\begin{eqnarray*}
\sum_i\length{\vx_i\vp} & = & \sum_i\phi_i(\vx_i-\vp) \\
& = & \sum_i\phi_i(\vx_i-\vq) + \sum_i\phi_i(\vq-\vp) \\
& \leq & \sum_i\length{\vx_i\vq} = \sum_i\length{\vx_i\vp}.
\end{eqnarray*}
Since there is equality throughout, it follows that each $\phi_i$ is a norming functional also of $\vx_i-\vq$.
Since $A$ is not collinear, we may choose an $\vx_i$ such that $\vx_i,\vp,\vq$ are not collinear.
Thus $\un{\vx_i-\vp}$ and $\un{\vx_i-\vq}$ are distinct unit vectors with the same norming functional $\phi_i$.
Thus $\un{\vx_i-\vp}\un{\vx_i-\vq}$ is a segment on the boundary of the unit ball, hence $X$ is not strictly convex.
\qed

We finally note the following difference between dimension two and higher dimensions.

\begin{theorem}[\cite{WH}]\label{3.21}
In any Minkowski plane $X$, for any finite $A\subset X$ we have $\conv A\cap\FT(A)\neq\emptyset$.
\end{theorem}

\begin{corollary}\label{sc2}
In a strictly convex Minkowski plane $X$, if $A\subset X$ is finite and non-collinear, the singleton $\FT(A)$ is always contained in $\conv A$.
\end{corollary}

\begin{theorem}[\cite{Durier, MR89f:90114}]\label{durier}
Let $\dim X\geq 3$.
Then the following are equivalent.
\begin{enumerate}
\item For any finite non-collinear $A\subset X$ we have $\FT(A)\subseteq\conv A$.\label{s1}
\item For any finite non-collinear $A\subset X$ we have $\FT(A)\cap\conv A\neq\emptyset$.\label{s2}
\item For any finite non-collinear $A\subset X$ we have $\FT(A)\subseteq\aff A$\label{s3}
\item For any finite non-collinear $A\subset X$ we have $\FT(A)\cap\aff A\neq\emptyset$.\label{s4}
\item $X$ is a Euclidean space.\label{s5}
\end{enumerate}
\end{theorem}
\proof
We sketch the proof as its non-trivial parts are in different references.
The implications \ref{s1}$\Rightarrow$\ref{s3}$\Rightarrow$\ref{s4} are trivial.
The well-known implication \ref{s5}$\Rightarrow$\ref{s1} is in \cite[Proposition~6.1]{KM}.
See Durier \cite{Durier} or Lewicki \cite{Lewicki} for \ref{s2}$\Rightarrow$\ref{s5}, as well as Ben\'\i tez, Fern\'andez and Soriano \cite{BFS} for stronger results.
Finally, see Cieslik \cite[Theorem~2.2]{MR89f:90114} for \ref{s4}$\Rightarrow$\ref{s2}.
\qed

\medskip
In the next section we characterize the situation in Minkowski planes when there are points of the FT locus that are outside the convex hull.

\section{Specific properties of FT loci in Minkowski planes}\label{specific}
In Euclidean space it is known that if $\FT(A)=\{\vp\}$, then $\vp\in A\cup\interior\conv A$ (see \cite[Proposition~6.1]{KM}.
In the light of Theorem~\ref{durier} we cannot hope for a similar statement in arbitrary Minkowski spaces of dimension at least three.
We now investigate to what extent we can have an analogue in Minkowski planes.
We first consider the case where $\FT(A)$ intersects the complement of $\conv A$.
\begin{definition}\label{4.1}
We say that a set $A=\{\vx_1,\dots,\vx_k,\vy_1,\dots,\vy_k\}$ forms a \emph{double cluster} with \emph{pairs} $\vx_i,\vy_i$ if $\un{\vx_i-\vy_i}$ are all contained in the same proper exposed face of the unit ball.
\end{definition}
Note that in the above definition, since we are in two dimensions, a proper exposed face of the unit ball is either a vertex (in which case the double cluster is necessarily a collinear set) or a segment.
Obviously, a double cluster forms $d$-concurrent $d$-segments.
By Corollary~\ref{dconc} the FT locus of a double cluster is a parallelogram with sides parallel to $\va$ and $\vb$, where $\va\vb$ is the exposed face of the unit ball in Definition~\ref{4.1}.

\begin{example}\label{doublecluster}
See Figure~\ref{fig5} for an example of a double cluster $A$ for which some FT points are not in $\conv A$.
\begin{figure}
\begin{center}
\includegraphics[scale = 0.7]{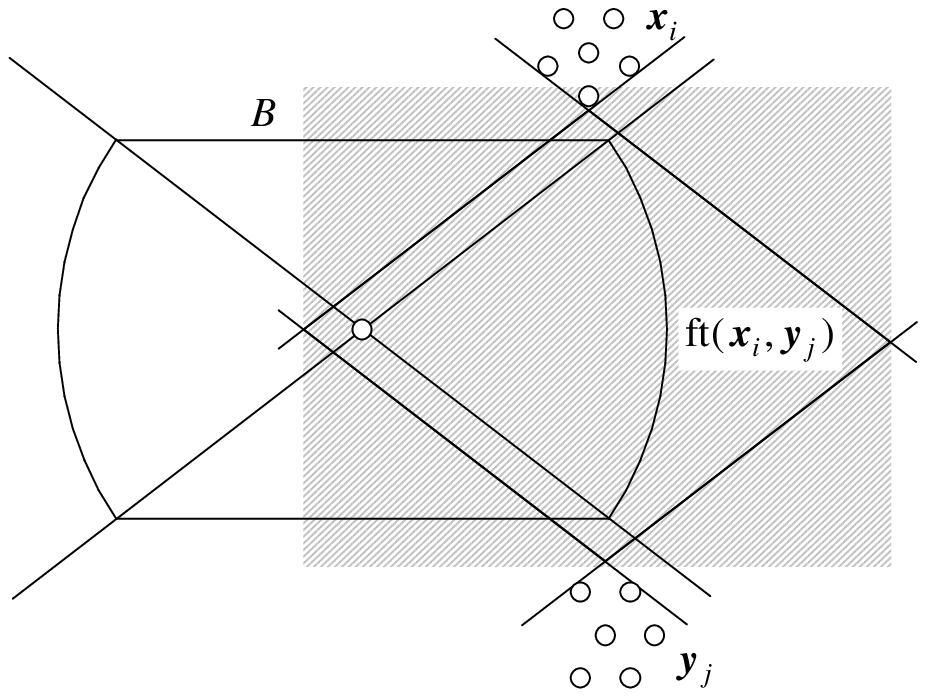}
\end{center}
\caption{A double cluster}\label{fig5}
\end{figure}
\end{example}

\begin{theorem}\label{dc}
Let $\vp$ be an FT point of a finite set $A$ in a Minkowski plane $X$ such that $\vp\not\in\conv A$.
Then $X$ is not strictly convex,
and $A$ is a double cluster.
In particular, $\FT(A)$ is the parallelogram $\bigcap_{i=1}^k\dsegment{\vx_i\vy_i}$.
\end{theorem}
\proof
Assume without loss of generality that $\vp=\vo$.
For each $\va_i\in A$, choose a norming functional $\phi_i$ such that $\sum_i\phi_i=o$.
Since $\vo\not\in\conv A$, we obtain that all $\phi_i$'s must be contained in a closed half plane bounded by a line $\ell$ through the origin in the dual.
Since the sum of the $\phi_i$'s is $o$, we must have that all the $\phi_i$'s must lie on $\ell$, and that there is an even number $2k$ of them, half being equal to some $\phi$, the other half to $-\phi$, say $\phi_1=\dots=\phi_k=\phi$, $\phi_{k+1}=\dots=\phi_{2k}=-\phi$.
Let $\vx_i=\va_i$ and $\vy_i=\va_{i+k}$ for all $i=1,\dots,k$.
It then follows that $\phi$ is a norming functional of any $\vx_i-\vy_j$, hence $\un{\vx_i-\vy_j}$ all lie on the same segment of the unit ball.
Since $\phi$ is a norming functional of all $\vx_i$ and $-\vy_j$, and $\vo\not\in\conv\{\vx_i,\vy_j\}$, we obtain that $X$ is not strictly convex.
\qed

As corollary we again obtain Corollary~\ref{sc2}, as well as
\begin{corollary}\label{4.4}
In a Minkowski plane, if $|A|$ is odd, then $\FT(A)\subseteq\conv A$.
In particular, neither a floating FT configuration of odd degree nor an absorbing FT configuration of even degree can be pointed.
\end{corollary}
A double cluster can give a pointed floating FT configuration of even degree, as in Example~\ref{doublecluster}.
However, by Corollary~\ref{sc2} this is not possible in strictly convex Minkowski planes.
Pointed \emph{absorbing} FT configurations of odd degree are possible, even in the Euclidean plane (see \cite[Remark~6.3.]{KM}).
By Corollary~\ref{4.4} this can only happen if the degree is odd in any Minkowski plane.
In the Euclidean plane one can say something stronger:
If $|A|$ is even and the FT point of $A$ is not in the interior of $\conv A$, then Proposition 6.2 of \cite{KM} gives that $A$ is \emph{almost collinear}, i.e.\ $A\setminus\{\vp\}$ is collinear for some $\vp\in A$.
We partially generalize this result to Minkowski planes.

\begin{definition}\label{4.6}
A set $A\subset X$ is a \emph{pseudo double cluster} if $A$ is the union of a double cluster $C$ together with an FT point of $C$ (called the \emph{centre} of the pseudo double cluster) and an arbitrary point.
\end{definition}
In a strictly convex Minkowski plane, a pseudo double cluster is an almost collinear set.

\begin{theorem}\label{prop62}
Let $A=\{\vx_0,\vx_1,\dots,\vx_{2k-1}\}$ in a Minkowski plane, and suppose that $\vx_0\in\FT(A)$ and $\vx_0$ is a vertex of $\conv A$.
Then $A$ is a pseudo double cluster with $\vx_0$ as centre.
\end{theorem}
\proof
Assume without loss of generality that $\vx_0=\vo$.
Choose norming functionals  $\phi_i$ of $\vx_i$, $i=1,\dots,2k-1$, such that $\norm{\sum\phi_i}\leq 1$.
Since $\vo$ is a vertex of $\conv A$, the $\phi_i$'s are all in a closed half plane bounded by a line $\ell$ through the origin in the dual.
Assume without loss of generality that the $\vx_i$'s are ordered such that the $\phi_i$'s are in order.
Let $\ell_1$ be a supporting line of the dual unit ball at $\phi_k$, and $\ell_0$ its parallel through $o$.
Let $H$ be the half plane bounded by $\ell_0$ containing $\phi_k$.
For any $i<k<j$, $\phi_i$ is on the closed arc from $-\phi_j$ to $\phi_k$ of the dual unit circle, otherwise $o$ is an interior point of the $\phi_i$'s, contradicting the fact that the $\phi_i$'s are in a closed half plane.
It follows that $\phi_i+\phi_j\in H$.
Thus $\sum_{i=1}^{2k-1}\phi_i\in\phi_k+H$.
Thus $\norm{\sum_{i=1}^{2k-1}\phi_i}\geq 1$, and therefore, $\norm{\sum_{i=1}^{2k-1}\phi_i}= 1$.
It follows that for all $i,j$ with $i<k<j$, $\phi_i+\phi_j\in\ell_0$.
Suppose $-\phi_j\neq\phi_i$.
Then $(-\phi_j)\phi_i$ is a segment on the boundary of the unit ball, parallel to $\ell_0$.
Thus $-\phi_j,\phi_i,\phi_k$ are collinear.
Similarly, $-\phi_i,\phi_j,\phi_k$ are collinear, a contradiction.

Thus $-\phi_j=\phi_i$ for all $i,j$ such that $i<k<j$.
Thus $\phi_1=\dots=\phi_{k-1}=\phi$ and $\phi_{k+1}=\dots=\phi_{2k-1}=-\phi$ for some unit $\phi\in X^*$.
It follows that the pairs $\vx_i, \vx_{i+k}$, $i=1,\dots,k-1$, form a double cluster, which has $o=\vx_0$ as an FT point.
\qed

\smallskip
The following is a complete generalization of \cite[Proposition~6.2]{KM} to strictly convex Minkowski planes.
\begin{corollary}\label{4.8}
Let $A=\{\vx_0,\vx_1,\dots,\vx_{2k-1}\}$ be given in a strictly convex Minkowski plane, and suppose that $\vx_0\in\FT(A)$ and $\vx_0$ is a vertex of $\conv A$.
Then for some $j=1,\dots,2k-1$ we have that $A\setminus\{\vx_j\}$ is a collinear set with $\vx_0$ as middle point.
\end{corollary}

The case where an FT point is on the edge, but is not a vertex of the convex hull of $A$, is more complicated.
What prevents the proof of Theorem~\ref{prop62} from going through in this case is that the unit vector parallel to such an edge may be a singular point of the unit ball.
If this cannot happen (such as when $X$ is smooth), then we have the following generalization of the above-mentioned \cite[Prop. 6.2]{KM}.
The proofs of the following two corollaries are simple adaptations of the proof of Theorem~\ref{prop62}.
\begin{corollary}\label{4.9}
In a smooth Minkowski plane, if $|A|$ is even and $\FT(A)$ intersects $\bd\conv A$, then $A$ is a double cluster or a pseudo double cluster.
\end{corollary}

\begin{corollary}\label{4.10}
In a smooth Minkowski plane, if $A$ is not a double cluster nor a pseudo double cluster, then $\FT(A)\subseteq A\cup\interior\conv A$.
\end{corollary}

If $X$ is not smooth, then it is always possible to find (even if $X$ is strictly convex) a set $A$ with an odd or an even number of points, which has an FT point on the interior of an edge of $\conv A$, as the following example shows.

\begin{example}\label{4.11}
Let $X$ be any non-smooth plane.
Let $\pm\vx_0$ be singular points on the boundary of the unit ball.
Let $\phi_0\phi_1$ be the set of all norming functionals of $\vp$.
Let $\vx_1$ be a unit vector with a norming functional parallel to $\phi_0\phi_1$.
Find unit functionals $\phi_2,\phi_3$ such that $\phi_3-\phi_2=\lambda(\phi_1-\phi_0)$ with $1\leq\lambda\leq 2$ and $\phi_3\neq-\phi_2$.
Let $\vx_2$ be a unit vector with $-\phi_2$ as norming functional, and $\vx_3$ a unit vector with $\phi_3$ as norming functional.
See Figure~\ref{fig6}.
\begin{figure}
\begin{center}
\includegraphics[scale = 0.85]{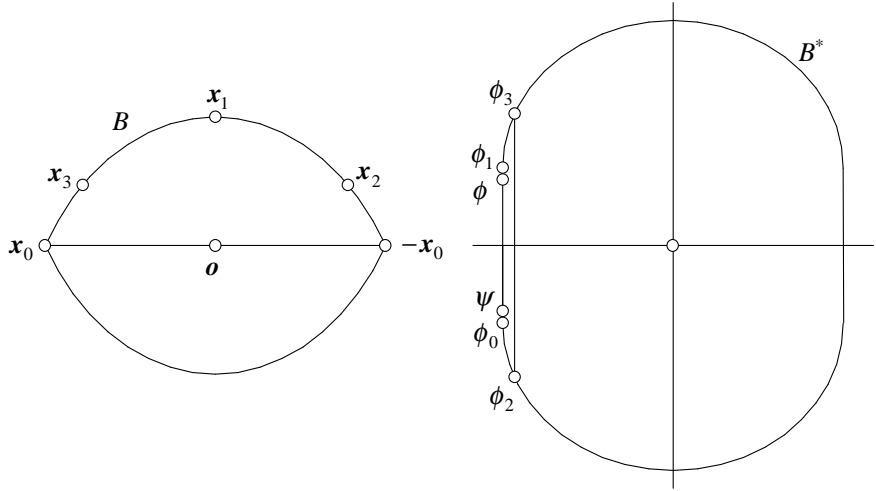}
\end{center}
\caption{An FT point on the interior of an edge of $\conv A$}\label{fig6}
\end{figure}
Choose $\phi,\psi\in\phi_1\phi_2$ such that $2(\phi-\psi)=\phi_3-\phi_2$.
Then $\vo$ is an FT point of e.g.\ $\{\pm\vx_0,\pm2\vx_0,\vx_2,\vx_3\}$.

It is also clear that we may add some odd number of multiples of $\vx_1$ and some further multiples of $\pm\vx_0$, and again obtain a set which has $\vo$ as FT point.
\end{example}

However, in the above example we still have an ``almost collinear'' situation in a weaker sense.
The following theorem shows that this necessarily happens.

\begin{theorem}\label{4.12}
In a Minkowski plane, if $A$ has an FT point disjoint from $A$, but on the interior of an edge of $\conv A$, then at least half of the points of $A$ must be on this edge.
\end{theorem}
\proof
Let $\vo$ be an FT point of $A$ between $\vx_0$ and $\vx_1$ in $A$.
Let $\phi_0\phi_1$ be the segment on the boundary of the dual unit ball containing all norming functionals of $\vx_0$.
Then all the norming functionals of points of $A$ that are not multiples of $\vx_0$ must be on the side of the line $\ell$ through $\phi_1$ parallel to $\phi_0+\phi_1$ opposite $\phi_0$ (or on $\ell$).
In order to obtain norming functionals of each $\vx\in A$ with sum $\vo$ we then must have at least as many multiples of $\vx_0$ as there are non-multiples.
\qed

Proposition 6.4 of \cite{KM} generalizes to all Minkowski planes as was shown in \cite{S5}:

\begin{theorem}[\cite{S5}]\label{4.13}
Let $\vp_0,\vp_1,\dots,\vp_n$ be distinct points in a Minkowski plane such that for any distinct $i,j$ satisfying $1\leq i,j\leq n$ the closed angle $\myangle \vp_i\vp_0\vp_j$ contains the reflection in $\vp_0$ of some $\vp_k$.
Then $n$ is necessarily odd and $\vp_0\in\FT(\vp_0,\dots,\vp_n)$.
\end{theorem}

For odd $n\leq 7$, and any convex $n$-gon $\vp_1\dots\vp_n$, there always exists $\vp_0$ such that the hypotheses of the above theorem is satisfied; see \cite{Tamvakis, BR}.
 
We now address the question of how many points of $A$ can be contained in $\FT(A)$.
Obviously, if $A$ consists of at most two points, then $A\subseteq\FT(A)$.
Examples~\ref{rectilinear} and \ref{reghex} show that it is possible for three and four points of $A$ to be in $\FT(A)$.

\begin{theorem}\label{twodim}
Let $X$ be a Minkowski plane and $A\subset X$.
Then $|A\cap\FT(A)|\leq 4$.
If $|A\cap\FT(A)|=4$, then $X$ has a parallelogram as unit ball, and $A$ contains a homothet of $\{\pm\vx\pm\vy\}$, where $\vx$ and $\vy$ are two consecutive vertices of the unit ball.
If $|A\cap\FT(A)|=3$, then $X$ has an affine regular hexagon as unit ball, and $A$ contains a homothet of $\{\vo,\vx,\vy\}$, where $\vx$ and $\vy$ are two consecutive vertices of the unit ball.
In both cases, $\FT(A)=\conv (A\cap\FT(A))$.
\end{theorem}

Before proving this theorem, we prove a technical lemma, and consider the higher-dimensional case.

\begin{lemma}\label{lemma2}
In any $d$-dimensional Minkowski space $X$, for each point $\vp\in A\cap\FT(A)$, $\vp$ is a vertex of $\conv (A\cap\FT(A))$, and $\{\un{\vq-\vp} : \vq\in A\cap\FT(A),\vq\neq \vp\}$ is contained in a proper exposed face of the unit ball.
\end{lemma}
\proof
Let $\vp\in A$.
By Theorem~\ref{DMchar} there exist norming functionals $\phi_\vx$ for each $\vx\in A\setminus\{\vp\}$ such that $\norm{\sum_{\vx\neq p}\phi_\vx}\leq 1$.
Then, for any $\vq\in\FT(A), \vq\neq \vp$, we have
\begin{eqnarray*}
\sum_{\vx\in A}\length{\vx\vq} & = & \sum_{\vx\in A,\vx\neq \vp}\length{\vx\vp} = \sum_{\vx\in A,\vx\neq \vp}\phi_\vx(\vx-\vp) \\
& = & \sum_{\vx\in A,\vx\neq \vp}\phi_\vx(\vx-\vq) + \sum_{\vx\in A,\vx\neq \vp}\phi_\vx(\vq-\vp)\\
& \leq & \sum_{\vx\in A,\vx\neq \vp}\length{\vx\vq} + \Bigl\|\sum_{\vx\in A,\vx\neq \vp}\phi_\vx\Bigr\|\length{\vq\vp}\\
& \leq & \sum_{\vx\in A,\vx\neq \vp}\length{\vx\vq} + \length{\vq\vp}\\
& = & \sum_{\vx\in A}\length{\vx\vq}.
\end{eqnarray*}
It follows that $\phi:=\sum_{\vx\in A,\vx\neq \vp}\phi_\vx$ is a norming functional of $\vq-\vp$.
Thus $\{\un{\vq-\vp}:\vq\in A\cap\FT(A),\vq\neq \vp\}$ is contained in the intersection of the unit ball with $\phi^{-1}(1)$.
This also means that $\phi$ strictly separates $\vp$ from $(A\cap\FT(A))\setminus\{\vp\}$, i.e.\ $\vp$ is a vertex of $\conv (A\cap\FT(A))$.
\qed

\begin{theorem}\label{highd}
Let $X$ be a $d$-dimensional Minkowski space, and $A\subset X$.
Then $|A\cap\FT(A)|\leq 2^d$.
If, furthermore, $|A\cap\FT(A)|=2^d$, then $X$ is isometric to a $d$-dimensional $L_1$ space, with $A\cap\FT(A)$ corresponding to a homothet of the Hamming cube $\{0,1\}^d$.
\end{theorem}
\proof
Let $C=A\cap\FT(A)$.
By Lemma~\ref{lemma2}, each point of $C$ is a vertex of $\conv C$.
We now show that $|C|\leq 2^d$.
For each $\vp\in C$, let 
\[C_\vp=\{\vx\in X:\vp+\lambda \vx\in\conv C\mbox{ for some }\lambda>0\}.\]
Each $C_\vp$ is a closed cone full-dimensional in the subspace $X'=\aff C - \aff C$, and also $C_\vp\cap-C_\vp=\{\vo\}$.
We now show that for any two distinct $\vp,\vq\in C$, $C_\vp\cap C_\vq$ does not have interior points (in $X'$) in common.
Note that $\vq-\vp\in C_\vp$ and $\vp-\vq\in C_\vq$.
Thus, if $C_\vp$ and $C_\vq$ have interior points in common, $\un{C_\vp}:=\{\un{\vx}:\vx\in C_\vp,\vx\neq \vo\}$ and $\un{C_\vq}$ are contained in the same proper exposed face of the unit ball.
But then $\un{\vp-\vq}$ and $\un{\vq-\vp}$ are both contained in this face, a contradiction.
It follows that $\{C_\vp:\vp\in C\}$ is a packing.
Thus $\{-\vp+\conv C:\vp\in C\}$ is a packing of translates of $\conv C$, all having the origin in common.
Let $B=\frac{1}{2}(\conv C-\conv C)$ be the central symmetrization of $\conv C$.
Then $\{B-\vp:\vp\in C\}$ is a packing of mutually touching translates of $B$ in $X'$.
By results of \cite{MR25:1488} it follows that $|C|\leq 2^d$, with equality iff $B$ is a $d$-dimensional parallelotope.
This happens only if $\conv C$ is a $d$-dimensional parallelotope, in which case $\bigcup_{\vp\in C}\un{C_\vp}$ forms the boundary of a cross-polytope.
\qed

\proof[Proof of Theorem~\ref{twodim}]
By Theorem~\ref{highd}, $|A\cap\FT(A)|\leq 4$, and if $|A\cap\FT(A)|=4$, then $X$ has a parallelogram as unit ball.

If $|A\cap\FT(A)|=3$, then $\conv (A\cap\FT(A))$ is a triangle.
The three cones in the proof of Theorem~\ref{highd} force the unit ball to be an affine regular hexagon.
\qed

\section{Centroids and Angles}\label{specific2}
We now consider characterizations of degree two absorbing and degree three floating FT configurations in Minkowski planes.
In studying such configurations, it is useful to introduce two special types of angles.
\begin{definition}\label{5.1}
An angle $\myangle \vx_1\vx_0\vx_2$ is \emph{critical} if there exists a point $\vx_3\neq \vx_0$ such that $\vx_0$ is an FT point of $\{\vx_1, \vx_2, \vx_3\}$.
\end{definition}
Critical angles are a direct generalization of Euclidean $120^\circ$ angles.
The ray $\Ray{\vx_0\vx_3}$ is unique for all critical angles iff the Minkowski plane is smooth and strictly convex \cite{DGGLW}.

The following characterizations of critical angles are well-known in the literature in the case of smooth, strictly convex planes \cite{CG, DGGLW}.
However, the generalization to arbitrary planes is simple.
\begin{lemma}\label{one}
The following are equivalent in a Minkowski plane.
\begin{enumerate}
\item $\myangle \vx_1\vx_0\vx_2$ is a critical angle.\label{onea}
\item If a circle \textup{(}in the norm\textup{)} with centre $\vx_0$ intersects the ray $\Ray{\vx_0\vx_1}$ at $\va$ and the ray $\Ray{\vx_0\vx_2}$ at $\vb$, then there exist lines $\ell_\va, \ell_\vb, \ell$ supporting the circle, $\ell_\va$ at $\va$ and $\ell_\vb$ at $\vb$, such that $\vx_0$ is the centroid of the triangle formed by $\ell_\va, \ell_\vb, \ell$.\label{oneb}
\item There exist norming functionals $\phi_1$ of $\vx_1-\vx_0$ and $\phi_2$ of $\vx_2-\vx_0$ such that $\norm{\phi_1+\phi_2} = 1$.\label{onec}
\end{enumerate}
\end{lemma}

\begin{definition}\label{5.2}
An angle $\myangle \vx_1\vx_0\vx_2$ is \emph{absorbing} if $\vx_0$ is an FT point of $\{\vx_0, \vx_1, \vx_2\}$.
\end{definition}
Thus $\{\vx_0\vx_1,\vx_0\vx_2\}$ is a degree two absorbing FT configuration iff $\myangle \vx_1\vx_0\vx_2$ is absorbing.
The following lemma furnishes a more direct description of absorbing angles.
In particular, an angle is absorbing iff it contains a critical angle.

\begin{lemma}\label{two}
The following are equivalent in a Minkowski plane.
\begin{enumerate}
\item $\myangle \vx_1\vx_0\vx_2$ is an absorbing angle.\label{twoa}
\item $\myangle \vx_1\vx_0\vx_2$ contains some critical angle $\myangle \vx_1'\vx_0\vx_2'$.\label{twob}
\item If a circle \textup{(}in the norm\textup{)} with centre $\vx_0$ intersects the ray $\Ray{\vx_0\vx_1}$ at $\va$ and the ray $\Ray{\vx_0\vx_2}$ at $\vb$, then there exist lines $\ell_\va, \ell_\vb, \ell$ with $\ell_\va$ supporting the circle at $\va$, $\ell_\vb$ at $\vb$, and $\ell$ not intersecting the interior of the circle, such that $\vx_0$ is the centroid of the triangle formed by $\ell_\va, \ell_\vb, \ell$.\label{twoc}
\item There exist norming functionals $\phi_1$ of $\vx_1-\vx_0$ and $\phi_2$ of $\vx_2-\vx_0$ such that $\norm{\phi_1+\phi_2} \leq 1$.\label{twod}
\end{enumerate}
\end{lemma}

\proof[Proof of Lemmas~\ref{one} and \ref{two}]
\ref{one}.\ref{onea}$\iff$\ref{one}.\ref{onec} and \ref{two}.\ref{twoa}$\iff$\ref{two}.\ref{twod} are straightforward using the subdifferential calculus.
\ref{one}.\ref{oneb}$\iff$\ref{one}.\ref{onec} and \ref{two}.\ref{twoc}$\iff$\ref{two}.\ref{twod} can be proved as in \cite{CG}.
\ref{two}.\ref{twoa}$\iff$\ref{two}.\ref{twob} follows from \ref{one}.\ref{onea}$\iff$\ref{one}.\ref{onec} and \ref{two}.\ref{twoa}$\iff$\ref{two}.\ref{twod}.
\qed

\smallskip
The following theorem, generalizing the characterization of degree three floating FT configurations in the Euclidean plane, is used in a characterization of the local structure of Steiner minimal trees \cite{S6}.
It is surprising that a Euclidean result can be completely generalized to all Minkowski planes.
It is again surprising that the proof is not simple.
Since we did not include a complete proof in \cite{S6}, we here prove the result in full.

\begin{theorem}\label{torithree}
The configuration $\{\vo\va_1, \vo\va_2, \vo\va_3\}$ is a floating FT configuration iff it is not pointed and all angles $\myangle \va_i\vo\va_j$ are critical.
\end{theorem}
\proof
$\Rightarrow$ The angles are all critical by definition.
By Corollary~\ref{4.4}, the configuration is not pointed.

$\Leftarrow$ Let $A_i$ be the set of norming functionals of $\va_i$ ($i=1,2,3$).
Note that if $\myangle\va_i\vo\va_j$ is a straight angle, then it is critical only if $A_i=-A_j$ is a non-degenerate segment.
Thus in all cases, $A_1,A_2,A_3$ are not contained in a closed half plane bounded by a line through the origin in the dual plane.
We now apply Theorem~\ref{DMchar} and the following lemma.
\qed

\begin{lemma}\label{threesegment}
In a Minkowski plane, let $A_i$ be the intersection of the unit ball with supporting line $\ell_i$ $(i=1,2,3)$ to the unit ball such that for all distinct $i,j\in\{1,2,3\}$ there exist $\va_i\in A_i$ and $\va_j\in A_j$ such that $\norm{\va_i+\va_j}=1$ and such that $A_1,A_2,A_3$ are not contained in a closed half plane bounded by a line through the origin.
Then there exist $\va_i\in A_i$, $(i=1,2,3)$ such that $\va_1+\va_2+\va_3=\vo$.
\end{lemma}
\proof
Assume without loss of generality that $\ell_1$ is vertical, $\ell_2$ horizontal, and $\ell_1, \ell_2$ enclose the unit square $[-1,1]^2\subseteq\R^2$.
Then the unit ball is contained in $[-1,1]^2$.
Let $A_1=\{1\}\times[-\alpha,\beta]$ and $A_2=[-\gamma,\delta]\times\{1\}$.
Then $A_1+A_2=[1-\gamma,1+\delta]\times[1-\alpha,1+\beta]$.
Note that $\alpha,\gamma\geq 0$: If $\alpha<0$ or $\gamma<0$, then $A_1+A_2$ would be outside $[-1,1]^2$, hence $\norm{\va_1+\va_2}>1$ for all $\va_1\in A_1$ and $\va_2\in A_2$.
For similar reasons, $A_3$ must contain a point with $x$-coordinate $\leq 0$, and a point with $y$-coordinate $\leq 0$.
Since $A_1+A_2$ intersects the unit ball, there must be a unit vector in $[1-\gamma,1]\times[1-\alpha,1]$.
Since $A_1,A_2,A_3$ are not in a closed half plane, one of the following cases must occur:
\begin{enumerate}
\item $-A_3$ intersects $[1-\gamma,1]\times[1-\alpha,1]$, in which case we are done.
\item $-A_3\subseteq [1-\gamma,1]\times[\beta,1-\alpha)$ (this is possible only if $\alpha+\beta<1$).
\item $-A_3\subseteq [\delta,1-\gamma)\times[1-\alpha,1]$ (this is possible only if $\gamma+\delta<1$).
\item $-A_3=A_1$.
\item $-A_3=A_2$.
\end{enumerate}
\emph{Case 2}.
\begin{eqnarray*}
A_1+A_3&\subseteq&\{1\}\times[-\alpha,\beta] + [-1,-1+\gamma]\times(-1+\alpha,\beta]\\
&=&[0,\gamma]\times(-1,0],
\end{eqnarray*}
which does not contain a unit vector, a contradiction.

\emph{Case 3}. A similar contradiction is obtained.

\emph{Case 4}. $A_1+A_3=\{0\}\times[-\alpha-\beta,\alpha+\beta]$.
Since $A_1+A_3$ must contain a unit vector, $1\in[-\alpha-\beta,\alpha+\beta]$, i.e.\ $\alpha+\beta\geq 1$.
Therefore $-\beta\leq -1+\alpha<\alpha$, and the vector $\opair{-1}{-1+\alpha}\in A_3$.
Then vectors $\opair{1}{-\alpha}\in A_1$, $\opair{0}{1}\in A_2$, $\opair{-1}{-1+\alpha}\in A_3$ have sum $\vo$.

\emph{Case 5}. Similar to Case 4.
\qed

\providecommand{\bysame}{\leavevmode\hbox to3em{\hrulefill}\thinspace}


\begin{thebibliography}{10}

\bibitem{BFS}
C.~Ben{\'\i}tez, M.~Fern{\'a}ndez, and M.~L. Soriano, \emph{Location of the
  {F}ermat centers of three points}, Submitted, 2000.

\bibitem{BMS}
V.~Boltyanski, H.~Martini, and P.~S. Soltan, \emph{Excursions into
  {C}ombinatorial {G}eometry}, Springer-Verlag, Berlin, 1997.

\bibitem{MR2000c:90002}
V.~Boltyanski, H.~Martini, and V.~Soltan, \emph{Geometric {M}ethods and
  {O}ptimization {P}roblems}, Kluwer Academic Publishers, Dordrecht, 1999.

\bibitem{BR}
M.~Bowron and S.~Rabinowitz, \emph{Solution to problem 10526}, Amer. Math.
  Monthly \textbf{104} (1997), 979--980.

\bibitem{CG}
G.~D. Chakerian and M.~A. Ghandehari, \emph{The {Fermat} problem in {Minkowski}
  spaces}, Geom. Dedicata \textbf{17} (1985), 227--238.

\bibitem{MR89f:90114}
D.~Cieslik, \emph{The {F}ermat-{S}teiner-{W}eber-problem in {M}inkowski
  spaces}, Optimization \textbf{19} (1988), 485--489.

\bibitem{Cieslik2}
D.~Cieslik, \emph{Steiner {M}inimal {T}rees}, Nonconvex {O}ptimization and its
  {A}pplications, vol.~23, Kluwer, Dordrecht, 1998.

\bibitem{MR25:1488}
L.~Danzer and B.~Gr{\"u}nbaum, \emph{\"{U}ber zwei {P}robleme bez\"uglich
  konvexer {K}\"orper von {P}. {E}rd{\H o}s und von {V}. {L}. {K}lee}, Math. Z.
  \textbf{79} (1962), 95--99.

\bibitem{MR1358610}
Z. Drezner (ed.), \emph{Facility {L}ocation, a {S}urvey of {A}pplications and
  {M}ethods}, Springer-Verlag, New York, 1995.

\bibitem{DGGLW}
D.-Z. Du, B.~Gao, R.~L. Graham, Z.-C. Liu, and P.-J. Wan, \emph{Minimum
  {Steiner} trees in normed planes}, Discrete Comput. Geom. \textbf{9} (1993),
  351--370.

\bibitem{Durier}
R.~Durier, \emph{The {F}ermat-{W}eber problem and inner product spaces}, J.
  Approx. Theory \textbf{78} (1994), no.~2, 161--173.

\bibitem{DM}
R.~Durier and C.~Michelot, \emph{Geometrical properties of the {Fermat}-{Weber}
  problem}, Europ. J. Oper. Res. \textbf{20} (1985), 332--343.

\bibitem{DM2}
\bysame, \emph{On the set of optimal points to the {W}eber problem},
  Transportation Sci. \textbf{28} (1994), 141--149.

\bibitem{KM}
Y.~S. Kupitz and H.~Martini, \emph{Geometric aspects of the generalized
  {Fermat}-{Torricelli} problem}, Intuitive Geometry, Bolyai Soc. Math. Stud.,
  vol.~6, 1997, pp.~55--127.

\bibitem{Lewicki}
G.~Lewicki, \emph{On a new proof of {D}urier's theorem}, Quaestiones Math.
  \textbf{18} (1995), 287--294.

\bibitem{Menger}
K.~Menger, \emph{Untersuchungen {\"u}ber allgemeine {M}etrik}, Math. Ann.
  \textbf{100} (1928), 75--163.

\bibitem{Rocka}
R.~T. Rockafellar, \emph{Convex {A}nalysis}, Princeton University Press,
  Princeton, 1997.

\bibitem{Schneider}
R. Schneider, \emph{Convex Bodies: the {Brunn}-{Minkowski} theory}, Encyclopedia of {Mathematics} and its {Applications} 44, Cambridge University Press, 1993.

\bibitem{S5}
K.~J. Swanepoel, \emph{Balancing unit vectors}, J. Combin. Theory Ser. A
  \textbf{89} (2000), 105--112.

\bibitem{S6}
\bysame, \emph{The local {S}teiner problem in normed planes}, Networks
  \textbf{36} (2000), 104--113.

\bibitem{Tamvakis}
H.~Tamvakis, \emph{Problem 10526}, Amer. Math. Monthly \textbf{103} (1996),
  427.

\bibitem{Thompson}
A.~C. Thompson, \emph{Minkowski {Geometry}}, Encyclopedia of {Mathematics} and
  its {Applications} 63, Cambridge University Press, 1996.

\bibitem{WH}
R.~E. Wendell and A.~P.~Hurter Jr., \emph{Location theory, dominance, and
  convexity}, Oper. Res. \textbf{21} (1973), 314--321.

\end{thebibliography}
\end{document}